\documentclass[journal,10pt]{IEEEtran}
\usepackage{ifpdf}

\ifCLASSINFOpdf
\else
\fi
\usepackage{amsmath}
\usepackage{algorithmicx}
\usepackage{algpseudocode}

\usepackage{array}

\usepackage{amsfonts}
\usepackage{amssymb}
\usepackage{graphicx}
\usepackage{float}
\usepackage{setspace}

\begin{document}
%
\title{Regularized maximum likelihood estimation of covariance matrices of elliptical distributions}


\author{\IEEEauthorblockN{Christophe Culan 
		,
		Claude Adnet\IEEEauthorrefmark{1}, 
	}
	\IEEEauthorblockA{ 
		Advanced Radar Concepts division,
		Thales, Limours 91338, France}
	\thanks{Manuscript received October XX, 2016; revised October XX, 2016. 
		Corresponding author: M. Culan (email: christophe.culan@thalesgroup.com).}}

\markboth{ArXiv prepublication}
{Culan \MakeLowercase{\textit{et al.}}: Regularized maximum likelihood estimators for elliptical distributions}
%



\IEEEtitleabstractindextext{%
\begin{abstract}
	The maximum likelihood principle is widely used in statistics, and the associated estimators often display good properties. indeed maximum likelihood estimators are guaranteed to be asymptotically efficient under mild conditions.\\
	However in some settings, one has too few samples to get a good estimation. It then becomes desirable to take into account prior information about the distribution one wants to estimate.\\
	One possible approach is to extend the maximum likelihood principle in a bayesian context, which then becomes a maximum a posteriori estimate; however this requires a distribution model on the distribution parameters.\\
	We shall therefore concentrate on the alternative approach of regularized estimators in this paper; we will show how they can be naturally introduced in the framework of maximum likelihood estimation, and how they can be extended to form robust estimators which can reject outliers.
\end{abstract}

}

\maketitle

\IEEEdisplaynontitleabstractindextext

%
\IEEEpeerreviewmaketitle

\section{Introduction: regularized estimation}

Many statistical estimators are shown to be asymptotically well behaved; in particular maximum likelihood estimators are consistent and efficient under relatively mild conditions, and often have good convergence properties, thus making them a popular choice for parametric estimation \cite{pratt1976fy}\cite{stigler2002statistics}\cite{hald1999history}\cite{aldrich1997ra}. However in a few sample regime, their performances can be greatly degraded, and in some cases one may not even be able to compute these estimators \cite{pfanzagl1994parametric}; this is notably the case for covariance estimations when the number of samples is of the same order or lower than the dimension of the space \cite{hoffbeck1996covariance}\cite{wiesel2012unified}.\\
In such situations, one has to add additional information in the estimation procedure, other than available samples. Hence a popular approach is to use regularized estimators, which are estimators purposely biased towards a prior distribution representing an initial knowledge of the model; the effect generally is a smoothing of the estimated parameter, although this largely depends upon the shape of the prior distribution; moreover this prior is not necessarily explicitly known. Various regularization methods have been proposed \cite{witten2009covariance}\cite{bickel2008covariance}\cite{barbaresco1996super}\cite{Duan16}\cite{Djafari02a}\cite{Ciuciu01}\cite{Giovannelli01b}. A popular approach is to derive such estimates by a penalized maximum likelihood approach \cite{xue2012positive}\cite{ravikumar2011high}\cite{wiesel2012unified}\cite{ollila2014regularized}.\\
The applied penalties are often chosen for convenience of computation; we introduce in this article a way to define such penalties which are in agreement with the theory of maximum likelihood.

\subsection{Notations and conventions}

In the following development, the following notations and conventions shall be observed:
\begin{itemize}
	\item For any topological space $\mathcal{R}_1$ and $\mathcal{R}_2$, the topological space $\mathcal{R}_1 \times \mathcal{R}_2$ is the product of $R_1$ and $R_2$, whereas $\mathcal{R}_1 \vee \mathcal{R}_2$ is the disjoint union of $R_1$ and $R_2$. One also notes, for a topological space $\mathcal{R}$,  $\mathcal{R}^N$ to be the product of $N$ copies of $\mathcal{R}$ and $\vee^N \mathcal{R}$ to be the disjoint union of $N$ copies of $\mathcal{R}$.
	\item Let $\mathcal{R}_1$, ..., $\mathcal{R}_N$ be topological spaces, and let $\mu_n$ be a measure defined on $\mathcal{R}_n$ for $1 \leq n \leq N$. The measure $\mu_1 \dots \mu_n = \prod_{n=1}^N \mu_n$ is the product measure of $(\mu_n)_{1\leq n\leq N}$ defined on $\prod_{n=1}^N \mathcal{R}_n$, whereas\\ $\mu_1 \vee \dots \vee \mu_n = \bigvee_{n=1}^N \mu_n$ is the measure of $\bigvee_{n=1}^N \mathcal{R}_n$ such that its restriction to $\mathcal{R}_n$ is $\mu_n$ for $1 \leq n \leq N$ and is the joint measure of $(\mu_n)_{1 \leq n \leq N}$.\\
	Moreover for a topological space $\mathcal{R}$ and a measure $\mu$ defined on $\mathcal{R}$, $\mu^N$ is the product measure of $N$ copies of $\mu$ defined on $\mathcal{R}^N$, whereas $\vee^N \mu$ is the joint measure of $N$ copies of $\mu$ defined on $N\cdot\mathcal{R}$.
	\item The vector space $\mathbb{C}^d$ is canonically identified to $\mathbb{R}^{2d}$ \cite{scharf1991statistical}.
	\item $\mathcal{H}_d(\mathbb{R})$ is the set of symmetric matrices of size $(d,d)$; $\mathcal{H}_d(\mathbb{C})$ is the set of hermitian matrices of size $(d,d)$.\\
	$\mathcal{H}_d^+(\mathbb{K})$ is the subset of matrices of 	$\mathcal{H}_d^+(\mathbb{K})$ which are positive; $\mathcal{HP}_d^+(\mathbb{K})$ is the subset of matrices of $\mathcal{H}_d^+(\mathbb{K})$ of unit determinant.
	\item The vector space $\mathbb{C}^d$ is canonically identified to $\mathbb{R}^{2d}$ \cite{scharf1991statistical}.
	\item  $X^\dagger$ is the transpose conjugate of any matrix or vector $X$.
	\item $\mu_{\mathbb{R}^d}$ is the Lebesgues measure of the space $\mathbb{R}^d$. Similarly, $\mu_{\mathbb{C}^d}$ is the Lebesgues measure of $\mathbb{C}^d \simeq \mathbb{R}^{2d}$.
	\item $\mathcal{S}_{d-1}$ is the $(d-1)$-sphere, which is identified as the following part of $\mathbb{R}^d$: $\mathcal{S}_{d-1} \simeq \left\{ x \in \mathbb{R}^d; x^\dagger x = 1 \right\}$. Similarly $\mathcal{S}_{2d-1}$ is identified to $\left\{ x \in \mathbb{C}^d; x^\dagger x = 1 \right\}$ in $\mathbb{C}^d$.
	\item $s_{d-1}$ shall denote the probability distribution of $\mathcal{S}_{d-1}$ isotropic for the canonical scalar product of $\mathbb{C}^d$, defined by:
	\begin{equation}
	\left<y,x\right> = \sum\limits_{k=0}^{d-1} \overline{y_k} x_k
	\end{equation}
	\item $\delta_x^\mathcal{R}$ is the Dirac delta distribution centered in $x$ in space $\mathcal{R}$.
	\item One shall note:
	\begin{equation}
	h_\delta = H\left(\delta_0^{\mathbb{R}}|\mu_{\mathbb{R}}\right) = \int_{t \in \mathbb{R}} \log\left(\frac{d\delta_0^{\mathbb{R}}}{dt}(t)\right)d\delta_0^{\mathbb{R}}(t)
	\end{equation}
	 It is the entropy of the Dirac distribution relative to the Lebesgues measure of $\mathbb{R}$ due to its translational symmetry, which is positive and infinite. It can be used to express different other entropies in higher dimensional settings:
	\begin{equation}
	\begin{array}{l}
	H\left(\delta_x^{\mathbb{R}^d}|\mu_{\mathbb{R}^d}\right) = d h_\delta\vspace{1.5 mm}\\
	H\left(\delta_x^{\mathbb{C}^d}|\mu_{\mathbb{C}^d}\right) = 2d h_\delta
	\end{array}
	\end{equation}
\end{itemize}

\section{Maximum likelihood with prior information}

\subsection{General theory}
\label{ssec:general_reg_ml}

Let us first recall that likelihood can be understood in terms of information; indeed it can be expressed as the relative entropy (or KL-divergence \cite{kullback1951information}) of the sampling distribution $S$ relative to the model distribution $P$ \cite{CVC_BT}:

\begin{equation}
l(P|S)= -H(S|P) = \int_{x \in \mathcal{R}} \log\left(\frac{dP}{dS}(x)\right)dS(x)
\end{equation}

Under the usual hypothesis in which one has $N$ i.i.d samples $(x_n)_{1\leq n \leq N}$, the usual sampling distribution is defined on a space $\vee^N \mathcal{B}$, with $\mathcal{B}$ being some base space, and is given by the disjoint union of Dirac distributions centered on each sample \cite{CVC_BT}:

\begin{equation}
S = \frac{1}{N}\bigvee_{n=1}^N \delta_{x_n}^\mathcal{B}
\end{equation}

The model distribution is given by $\frac{1}{N} \vee^N P_\mathcal{B}$, with $P_\mathcal{B}$ being some distribution of $\mathcal{B}$.\\

Thus the likelihood can be expressed under this hypothesis by:

\begin{equation}
l(P|S) = -\frac{1}{N} \sum_{n=1}^N H(\delta_{x_n}^\mathcal{B}|P_\mathcal{B})
\end{equation}

Obviously this particular instantiation of the maximum likelihood problem is not well suited for our current problem, as we would like to take into account some prior information in addition to the information provided by samples.\\ 
We shall focus on a situation in which we have an existing prior estimate, and we want to take into account this prior in our estimation procedure \cite{wiesel2012unified}\cite{ollila2014regularized}.\\
Fortunately the entropic formulation of the maximum likelihood is quite flexible; in particular one can use different models for the sampling distribution, as is done for example in \cite{CVC_BT} to take into account the circular symmetry of the distributions.\\
We therefore propose to take into account the prior knowledge as part of the sampling distribution. The sampling distribution under the i.i.d hypothesis, with additional prior information can then be modeled as:\\

\begin{equation}
S = \left(\frac{\alpha}{N}\bigvee_{n=1}^N S_n\right) \vee \left((1-\alpha)P_\text{prior}\right)  
\end{equation}

with:

\begin{itemize}
	\item $S_n$ the sample distribution associated to the n-th sample, defined on a base space $\mathcal{B}$; it would generally be given by $\delta_{x_n}^\mathcal{B}$.
	\item $P_\text{prior}$ the prior distribution. It is not necessarily defined on the same space as the sample distributions depending on the information it carries, and might be defined on a quotiented image of the base space, which shall be noted as $[\mathcal{B}]$.
	\item $\alpha \in [0;1]$ the integration factor, which represents how much ones trusts the prior estimate compared to the newly acquired samples.
\end{itemize}

The distribution model one wants to fit can on the other hand be expressed as:

\begin{equation}
\left(\frac{\alpha}{N}\vee^N P\right)\vee \left((1-\alpha)[P]\right)
\end{equation}

with:

\begin{itemize}
	\item $P$ the distribution model
	\item $[P]$ the distribution model resulting from $P$ on the quotiented space $[\mathcal{B}]$
\end{itemize}

The corresponding likelihood can therefore be expressed as:

\begin{equation}
l(P|S) = -\frac{\alpha}{N}\sum\limits_{n=1}^N H(S_n|P) - (1-\alpha)H(P_\text{prior}|[P])
\end{equation}

Therefore this results in a penalized version of maximum likelihood , the penalty being defined from the relative entropy of the prior distribution for the model distribution \cite{ollila2014regularized}\cite{wiesel2012unified}.

\section{Regularized covariance estimators}

We shall now apply the regularized maximum likelihood equation to different covariance estimation problems. We'll focus on the two following cases:

\begin{itemize}
	\item The distribution which we try to estimate is elliptically symmetric (ES/CES) \cite{chmielewski1981elliptically}\cite{ollila2012complex}. No hypothesis is made on the radial distribution of this distribution; this corresponds to the conditions for which Tyler's estimator is the maximum likelihood estimator for i.i.d samples \cite{CVC_BT}\cite{tyler1987statistical}\cite{pascal2008covariance}.
	\item The distribution is gaussian, or circular gaussian for complex distributions. Note that the corresponding maximum likelihood estimators under the i.i.d hypothesis are given in \cite{CVC_BT}.
\end{itemize}

\subsection{Regularized Tyler's estimator}

Let us first recall that an elliptical distribution can be expressed in general by \cite{CVC_BT}:

\begin{equation}
dP(x) = \frac{1}{\sqrt{\left|R\right|}}dQ\left(\sqrt{x^\dagger R^{-1}x}\right)ds_{d-1}\left( \frac{L(R)^{-1} x}{\sqrt{x^\dagger R^{-1} x}} \right)
\end{equation}

with the radial distribution $Q$ being any probability measure on $\mathbb{R}_+$, and the correlation matrix $R$ being any matrix of $\mathcal{HP}_d^+(\mathbb{K})$ \cite{CVC_BT}.

Supposing here that one has no a priori on the radial distribution, the prior distribution one wants to use should be defined on the unit sphere, and corresponds to the unique probability distribution of the unit sphere which is isotropic for the corresponding correlation matrix $R_\text{prior}$, noted as $s_{d-1}^{R_\text{prior}}$. The entropy between two such distributions isotropic for correlation matrices $R_1$ and $R_2$ shall be noted $H(R_1|R_2)$. Thus the concentrated likelihood for samples $\left(x_n\right)_{1\leq n \leq N}$in $\mathbb{R}^d$ under i.i.d hypothesis is given, after maximization over the radial distribution, by \cite{CVC_BT}:

\begin{equation}
\begin{split}
l(R|S) =& -dh_\delta-\frac{\alpha(d-1)}{2N}\sum\limits_{n=1}^N \log\left({x_n}^\dagger R^{-1} x_n\right) \\
&-(1-\alpha)H(R_\text{prior}|R)
\end{split}
\end{equation}

$H(R_\text{prior}|R)$ however represents a challenge to compute. Fortunately, one can characterize its differential with respect to $R$ (see \ref{sec:KL_angular} for more details):\\

\begin{flalign*}
dH(R_\text{prior}|R) &&
\end{flalign*}
\begin{equation}
= \frac{d-1}{2}\text{tr}\left(\left(R^{\frac{1}{2}} \text{anscm}(R^{-\frac{1}{2}}R_\text{prior}R^{-\frac{1}{2}})R^{\frac{1}{2}}\right)d\left(R^{-1}\right)\right)
\end{equation}

with $\text{anscm}(\Sigma)$ being the expected value of the normalized sample covariance estimator (NSCM) for an elliptical distribution of correlation matrix proportional to $\Sigma$.\\
Unfortunately there is no explicit way of computing this matrix in the general case to the authors' knowledge and one has to resort to numerical methods to extract its eigenvalues. However it has been characterized for complex circularly symmetric distribution in \cite{bausson2007first}.

In any case this leads to the following maximum likelihood equation:

\begin{equation}
\begin{split}
\text{tr}\bigg(\bigg(\lambda R\\
& -\frac{(1-\alpha)(d-1)}{2}R^{\frac{1}{2}} \text{anscm}\left(R^{-\frac{1}{2}} R_\text{prior} R^{-\frac{1}{2}}\right)R^{\frac{1}{2}}\\
&-\frac{\alpha(d-1)}{2N}\sum\limits_{n=1}^N \frac{x_n{x_n}^\dagger}{{x_n}R^{-1}x_n} \bigg)d\left(R^{-1}\right)\bigg)=0
\end{split}
\end{equation}

Without any structural constrain on $R$, this simplifies to:

\begin{equation}
\begin{split}
R =& \frac{1}{2\lambda(d-1)} \bigg( (1-\alpha)R^\frac{1}{2} \text{anscm}\left(R^{-\frac{1}{2}} R_\text{prior} R^{-\frac{1}{2}}\right) R^\frac{1}{2}\\
&+\frac{\alpha}{N}\sum\limits_{n=1}^N \frac{x_n {x_n}^\dagger}{{x_n}^\dagger R^{-1} x_n} \bigg)
\end{split}
\end{equation}

This can be solved by using the following numerical procedure:\\

\begin{onehalfspace}
	\begin{samepage}
		\begin{algorithmic}[1]
			\Function{reg\_Tyler}{$R_\text{prior}$,$\left(x_n\right)_{1 \leq n \leq N}$,$\alpha$,$\epsilon$}
			\State $R \gets R_\text{prior}$
			\Repeat
			\State $R_\frac{1}{2} \gets \sqrt{R}$
			\State $R_{-\frac{1}{2}} \gets {R_\frac{1}{2}}^{-1}$\vspace{1.5 mm}
			\State $S \gets (1-\alpha)\text{anscm}\left(R_{-\frac{1}{2}}R_\text{prior}R_{-\frac{1}{2}}\right)$
			\begin{align*}
			&&+\frac{\alpha }{N}R_{-\frac{1}{2}}\sum\limits_{n=1}^N \frac{x_n {x_n}^\dagger}{{x_n}^\dagger R_{-\frac{1}{2}}R_{-\frac{1}{2}} x_n}R_{-\frac{1}{2}}
			\end{align*}
			\State $S \gets \frac{d}{\text{tr}(S)}S$\vspace{1.5 mm}
			\State $R \gets R_\frac{1}{2}\exp(S-I)R_\frac{1}{2}$\vspace{1.5 mm}
			\State $R \gets \frac{R}{\text{tr}(R)}$\vspace{1.5 mm}
			\Until{$\text{tr}\left(\left(R_{-1} R - I\right)^2\right) \leq \epsilon$}\vspace{1.5 mm}
			\State \Return $R_{-\frac{1}{2}}$
			\EndFunction
		\end{algorithmic}
	\end{samepage}
\end{onehalfspace}
\vspace{2 mm}

In the case of a distribution in $\mathbb{C}^d$ under the circularity hypothesis, the maximum likelihood equation under i.i.d hypothesis is slightly different; indeed the sampling distribution is then given by:

\begin{equation}
S = \frac{\alpha}{N} \bigvee_{n=1}^N \delta_{[x_n]}^\mathbb{C^d} + (1-\alpha)s_{d-1}^{R_\text{prior}}
\end{equation}

By deriving the concentrated likelihood in a manner similar to \cite{CVC_BT} and differentiating it, on gets the following equation:

\begin{equation}
\begin{split}
\text{tr}\bigg(\bigg(\lambda R\\
&-(1-\alpha)\left(d-\frac{1}{2}\right)R^\frac{1}{2} \text{anscm}\left(R^{-\frac{1}{2}}R_\text{prior}R^{-\frac{1}{2}}\right)R^\frac{1}{2}\\
&-\frac{\alpha(d-1)}{N}\sum\limits_{n=1}^N \frac{x_n{x_n}^\dagger}{{x_n}^\dagger R^{-1} x_n}
\bigg)d\left(R^{-1}\right)\bigg)	= 0
\end{split}
\end{equation}		

This is almost the same equation, although each sample actually carries less information compared to the prior distribution for a same $\alpha$ parameter.\\
The equation can thus be solved numerically by calling the same algorithm \textproc{reg\_Tyler} with a different $\alpha$ parameter:\\ 

\textproc{reg\_Tyler}($R_\text{prior}$,$\left(x_n\right)_{1 \leq n \leq N}$,$\frac{\alpha(d-1)}{d-\frac{1}{2}(1+\alpha)}$,$\epsilon$)\\

Note however that the computation of the ansm function differs in the real and complex circular cases (see annex).\\

These regularized Tyler estimators are more complicated to compute than those found in the existing litterature by applying different penalties (specifically ones maintaining the convexity of the problem) and offer no theoretical guaranty of convergence \cite{wiesel2012geodesic}\cite{ollila2014regularized}\cite{sun2014regularized}. However they use a penalty which is consistent with the definition of the likelihood as an information theoretic quantity, which is the aim of this paper.

\subsection{Regularized estimators for gaussian distributions}
\label{ssec:reg_gauss}
A subclass of interest of elliptical models are estimations constrained on a subclass of models corresponding to a fixed scaled radial distribution, of the from:

\begin{equation}
dP(x) = \sigma dQ_0(\sigma r)ds_{d-1}\left( \frac{L(R)^{-1} x}{\sqrt{x^\dagger R^{-1} x}} \right)
\end{equation}

Indeed most M-estimators can be expressed as maximum likelihood estimators for one of these subclasses of models \cite{CVC_BT}.\\
However the existence of an easy to compute regularized estimator then depends on the expression of the gradient of $H(\Sigma_1|\Sigma_2)$ with respect to $\Sigma_2$ for two distributions with the same base radial distribution $Q_0$. Unfortunately there is no analytical expression of this gradient except in specific cases \cite{hosseini2016inference}. We shall therefore focus on one of these specific cases: that of gaussian distributions, which corresponds to $Q_0$ following a $\chi^2(c_\mathbb{K}d)$ distribution, with $c_\mathbb{K} = 1$ for $\mathbb{K} = \mathbb{R}$, and $c_\mathbb{K} = 2$ for $\mathbb{K} = \mathbb{C}$.\\

The relative entropy of two such gaussian models is given by \cite{duchi2007derivations}:

\begin{equation}
H(\Sigma_1|\Sigma_2) = \frac{c_\mathbb{K}}{2}\left( \text{tr}\left({\Sigma_2}^{-1}\Sigma_1\right) + \log|\Sigma_2|-\log|\Sigma_1|-d\right)
\end{equation}

with $\Sigma_1$ and $\Sigma_2$ being the covariance matrices of said models.
In the standard case of gaussian distributions in $\mathbb{R}^d$, the resulting likelihood is thus given up to constant additive terms by:

\begin{equation}
\begin{split}
l(\Sigma|S) =& -\frac{\alpha}{2N}\sum\limits_{n=1}^N {x_n}^\dagger \Sigma^{-1} x_n \\
&-\frac{1-\alpha}{2}\text{tr}\left(\Sigma^{-1}\Sigma_\text{prior}\right)-\frac{1}{2}\log|\Sigma|	
\end{split}
\end{equation}

The corresponding maximum likelihood estimator is then given by a weighted mean of the covariance of the prior distribution and that of the sample covariance matrix of the new samples (SCM):

\begin{equation}
\Sigma = (1-\alpha)\Sigma_\text{prior}+\frac{\alpha}{N}\sum\limits_{n=1}^N x_n {x_n}^\dagger
\end{equation}

In the complex circular case, one can take into account the phase symmetry of the sampling distribution, in a manner similar to Tyler's estimator; this results in the following likelihood, up to constant additive terms:

\begin{equation}
\begin{split}
l(\Sigma|S) = &-\frac{\alpha}{N}\sum\limits_{n=1}^N \left({x_n}^\dagger \Sigma^{-1} x_n -\frac{1}{2}\log\left({x_n}^\dagger \Sigma^{-1} x_n\right)\right)\\
&-(1-\alpha)\text{tr}\left(\Sigma^{-1}\Sigma_\text{prior}\right)-\left(1-\frac{\alpha}{2d}\right)\log|\Sigma|
\end{split}
\end{equation}

This can be differentiated, leading to the following maximum likelihood equation:

\begin{equation}
\begin{split}
\Sigma = &\frac{1}{1-\frac{\alpha}{2d}}( (1-\alpha)\Sigma_\text{prior}\\
&+\frac{\alpha}{N}\sum\limits_{n=1}^N \left(1-\frac{1}{2{x_n}^\dagger \Sigma^{-1} x_n}\right)x_n{x_n}^\dagger )
\end{split}
\end{equation}

A solution can be found numerically by adapting the \textproc{cg\_cov} algorithm described in \cite{CVC_BT}:\\

\begin{onehalfspace}
	\begin{samepage}
		\begin{algorithmic}[1]
			\Function{reg\_cg\_cov}{$\Sigma_\text{prior},(x_n)_{1 \leq n \leq N},\alpha,\epsilon,K_{max}$}\vspace{1.5 mm}
			\State $\Sigma \gets (1-\alpha)\Sigma_\text{prior}+\frac{\alpha}{N}\sum\limits_{n=1}^{N} x_n {x_n}^\dagger$\vspace{1.5 mm}
			\Repeat
			\State $\Sigma_{\frac{1}{2}} \gets \sqrt{\Sigma}$\vspace{1.5 mm}
			\State $\Sigma_{-\frac{1}{2}} \gets {\Sigma_{\frac{1}{2}}}^{-1}$\vspace{1.5 mm}
			\State $S \gets \frac{1}{\left(1-\frac{\alpha}{2d}\right)}( (1-\alpha)\Sigma_\text{prior}$
			\begin{align*}
			&&+\frac{\alpha}{N}\sum\limits_{n=1}^N \left(1-\frac{1}{2{x_n}^\dagger {\Sigma_{-1}} x_n}\right)x_n{x_n}^\dagger )
			\end{align*}
			\State $R \gets \Sigma_{\frac{1}{2}}\exp\left(\Sigma_{-\frac{1}{2}} S \Sigma_{-\frac{1}{2}}-I\right)\Sigma_{\frac{1}{2}}$\vspace{1.5 mm}
			\Until{$\text{tr}\left(\left(\Sigma_{-\frac{1}{2}} \Sigma \Sigma_{-\frac{1}{2}} - I\right)^2\right) \leq \epsilon$}\vspace{1.5 mm}
			
			\State \Return $\Sigma_{-\frac{1}{2}}$
			\EndFunction
		\end{algorithmic}
	\end{samepage}
\end{onehalfspace}
\vspace{2 mm}

Thus in the specific case of gaussian and complex circular gaussian variables, the regularized maximum likelihood estimators can be computed efficiently.

\section{Partial estimation and regularized estimators}

We have thus successfully extended maximum likelihood estimators of covariance matrices to a regularized estimation setting. We shall now discuss a further extension of these regularized estimators in order to improve their robustness to cases of contamination by outliers. In particular, we would like to extend the partial maximum likelihood estimation method developed in \cite{CVC_partial} to these estimators.\\

Let us recall that partial estimation as introduced in \cite{CVC_partial} consists in the maximization of a partial likelihood:

\begin{equation}
l_{\left| X\right.}(P|S) = -H_{\left| X\right.}(P|S) = \int_{x \in X} \log\left(\frac{dP}{dS}(x)\right)dS(x)
\end{equation}

under the constrain that $X$ can be any measurable set such that $S(X) \geq p$ with  $p\in \left[ 0; 1\right]$ being the order of the partial estimation.\\

Unfortunately this definition is impractical in the case of regularized estimators for which $\mathcal{R} = \vee^N \mathcal{B} \vee [\mathcal{B}]$. Indeed one would have to consider partial entropy integrals between the prior and model distribution of the form:

\begin{equation}
	H_{\left| X\cap [\mathcal{B}]\right.}(P_\text{prior}|P) = \int_{x \in X \cap[\mathcal{B}]} \log\left(\frac{dP_\text{prior}}{dP}(x)\right)dP_\text{prior}(x)
\end{equation}

We therefore propose to restrict the possible partial domains to subsets $X$ of the form $X \vee [\mathcal{B}]$, with $X \subset \vee^N \mathcal{B}$ and such that $S(X) \geq (1-\alpha) + \alpha p$. Thus the corresponding restricted partial likelihood, maximized on the partial domain $X$ is now given by \cite{CVC_partial}: 

\begin{equation}
\begin{split}
l_p(P|S) =& -H_p(P|S)\\ 
=& \sup\limits_{S(X) \geq \alpha p+(1-\alpha)} \int_{x \in X} \log\left(\frac{dP}{dS}(x)\right)dS(x)
\end{split}
\end{equation}

This is the restricted partial likelihood of order $p$, conditional to the fact that the prior distribution $P_\text{prior}$ contains no outliers.

This can be recast as:

\begin{equation}
\begin{split}
l_p(P|S) =& -\frac{\alpha}{\lceil pN \rceil}\sum\limits_{n=1}^{\lceil pN \rceil} H\left(\delta_{x_{o(n)}}^\mathcal{B}|P\right)\\
&-(1-\alpha)H(P_\text{prior}|[P])
\end{split}
\end{equation}

This is a penalized form of the partial likelihood as expressed in \cite{CVC_partial}. Thus the partial estimation procedure introduced in \cite{CVC_partial} can be directly extended to such regularized estimators.

\subsection{Restricted partial maximum likelihood principle for regularized estimators}

We shall now apply the partial estimation procedure to the different regularized estimators introduced so far.

Outlined below is the algorithmic procedures corresponding to the partial version of the regularized Tyler estimator for a distribution in $\mathbb{R}^d$:\\

\begin{onehalfspace}
	\begin{samepage}
		\begin{algorithmic}[1]
			\Function{reg\_pTyler}{$R_\text{prior},(x_n)_{1 \leq n \leq N},\alpha,p,\epsilon,K_\text{max}$}
			\State $R_\frac{1}{2} \gets \sqrt{R_\text{prior}}$
			\State $R_{-\frac{1}{2}} \gets {R_\frac{1}{2}}^{-1}$ 
			\For{$k$ \textbf{from} $1$ \textbf{to} $K_\text{max}$}
			\For{$n$ \textbf{from} $1$ \textbf{to} $N$}
			\State $\tau_n \gets {x_n}^\dagger R_{-\frac{1}{2}}R_{-\frac{1}{2}} x_n$
			\EndFor\vspace{1.5 mm}
			\State $o \gets \text{argsort}_{\uparrow}\left(\left(\tau_n\right)_{1 \leq n \leq N}\right)$\vspace{1.5 mm} 
			\State $S \gets (1-\alpha)\text{anscm}\left(R_{-\frac{1}{2}}R_\text{prior}R_{-\frac{1}{2}}\right)$
			\begin{align*}
			&&+\frac{\alpha}{\lceil pN \rceil}R_{-\frac{1}{2}}\sum\limits_{n=1}^{\lceil pN \rceil} \frac{x_{o(n)} {x_{o(n)}}^\dagger}{\tau_{o(n)}}R_{-\frac{1}{2}}
			\end{align*}
			\State $S \gets \frac{d}{\text{tr}(S)}S$\vspace{1.5 mm}
			\State $R \gets R_\frac{1}{2}\exp(S-I)R_\frac{1}{2}$\vspace{1.5 mm}
			\State $R \gets \frac{R}{\text{tr}(R)}$\vspace{1.5 mm}
			\If{$\text{tr}\left(\left(R_{-1} R - I\right)^2\right) \leq \epsilon$}\vspace{1.5 mm}
			\State \textbf{break}	
			\Else	
			\State $R_\frac{1}{2} \gets \sqrt{R}$
			\State $R_{-\frac{1}{2}} \gets {R_\frac{1}{2}}^{-1}$
			\EndIf
			\EndFor
			
			\item \textbf{return} $R_{-\frac{1}{2}}$
			\EndFunction
		\end{algorithmic}
	\end{samepage}
\end{onehalfspace}
\vspace{2 mm}

The corresponding partial procedure for a distribution in $\mathbb{C}^d$ is the same as Tyler but using a different $\alpha$ parameter:\\

\textproc{reg\_pTyler}($R_\text{prior}$,$\left(x_n\right)_{1 \leq n \leq N}$,$\frac{\alpha(d-1)}{d-\frac{1}{2}(1+\alpha)}$,$p$,$\epsilon$)\\

As stated in \ref{ssec:reg_gauss}, it is in general quite difficult to derive regularized estimators for M-estimators which are consistent with maximum likelihood theory as stated in \ref{ssec:general_reg_ml}except for some specific cases. The same remains valid for partial regularized estimates, for which we shall again concentrate on the specific case of gaussian models. The corresponding partial regularized estimators, corresponding to the partial versions of the regularized sample covariance in the real case and the \textproc{reg\_cg\_cov} algorithm in the complex circular case are outlined below:\\

\begin{onehalfspace}
	\begin{samepage}
		\begin{algorithmic}[1]
			\Function{reg\_partial\_SCM}{$\Sigma_\text{prior},\left(x_n\right)_{1 \leq n \leq N},\alpha,K_\text{max}$}
			\State $\Sigma \gets \Sigma_\text{prior}$
			\For{$n$ \textbf{from} $1$ \textbf{to} $N$}
			\State $\tau_n \gets {x_n}^\dagger \Sigma^{-1} x_n$
			\EndFor\vspace{1.5 mm}
			\State $o_0 \gets \text{argsort}_{\uparrow}\left(\left(x_n\right)_{1 \leq n \leq N}\right)$\vspace{1.5 mm}
			\For{$k$ \textbf{from} $1$ \textbf{to} $K_\text{max}$}\vspace{1.5 mm}
			\State $\Sigma \gets (1-\alpha)\Sigma_\text{prior}
			+\frac{\alpha}{\lceil pN \rceil}\sum\limits_{n=1}^{\lceil pN \rceil} x_{o_0(n)} {x_{o_0(n)}}^\dagger$\vspace{1.5 mm}
			\For{$n$ \textbf{from} $1$ \textbf{to} $N$}
			\State $\tau_n \gets {x_n}^\dagger \Sigma^{-1} x_n$
			\EndFor\vspace{1.5 mm}
			\State $o_1 \gets \text{argsort}_{\uparrow}\left(\left(\tau_n\right)_{1 \leq n \leq N}\right)$\vspace{1.5 mm}
			\If{$o_0 = o_1$}
			\State \textbf{break}
			\Else
			\State $o_0 \gets o_1$
			\EndIf
			\EndFor
			\State \Return $\Sigma$
			\EndFunction
		\end{algorithmic}
	\end{samepage}
\end{onehalfspace}
\vspace{2 mm}

In the complex circular case the corresponding algorithm is given by:\\

\begin{onehalfspace}
	\begin{samepage}
		\begin{algorithmic}[1]
			\Function{reg\_pcg\_cov}{$\Sigma_\text{prior}\left(x_n\right)_{1 \leq n \leq N},\alpha,\epsilon,K_\text{max}$}\vspace{1.5 mm}
			\State $\Sigma \gets (1-\alpha) \Sigma_\text{prior}+\frac{\alpha}{N}\sum\limits_{n=1}^{N} x_n {x_n}^\dagger$\vspace{1.5 mm}
			\State $\Sigma_{\frac{1}{2}} \gets \sqrt{\Sigma}$
			\State $\Sigma_{-\frac{1}{2}} \gets {\Sigma_{\frac{1}{2}}}^{-1}$
			\For{$k$ \textbf{from} $1$ \textbf{to} $K_\text{max}$}
			\For{$n$ \textbf{from} $1$ \textbf{to} $N$}
			\State $\tau_n \gets {x_n}^\dagger \Sigma^{-1} x_n$
			\EndFor\vspace{1.5 mm}
			\State $o \gets \text{argsort}_{\uparrow}\left(\left(\tau_n-\frac{1}{2}\log\left(\tau_n\right)\right)_{1 \leq n \leq N}\right)$\vspace{1.5 mm}
			\State $S \gets \frac{1}{\left(1-\frac{\alpha}{2d}\right)} \bigg( (1-\alpha)\Sigma_\text{prior}$
			\begin{align*}
			&&+\frac{\alpha}{\lceil pN \rceil}\sum\limits_{n=1}^{\lceil pN \rceil} \left(1-\frac{1}{2\tau_{o(n)}}\right)x_{o(n)}{x_{o(n)}}^\dagger \bigg)
			\end{align*}
			\State $R \gets \Sigma_{\frac{1}{2}}\exp\left(\Sigma_{-\frac{1}{2}} S \Sigma_{-\frac{1}{2}}\right)\Sigma_{\frac{1}{2}}$\vspace{1.5 mm}
			\If{$\text{tr}\left(\left(\Sigma_{-\frac{1}{2}} \Sigma \Sigma_{-\frac{1}{2}} - I\right)^2\right) \leq \epsilon$}\vspace{1.5 mm} 
			\State \textbf{break}	
			\Else
			\State $\Sigma_{\frac{1}{2}} \gets \sqrt{\Sigma}$
			\State $\Sigma_{-\frac{1}{2}} \gets {\Sigma_{\frac{1}{2}}}^{-1}$
			\EndIf
			\EndFor
			\State \Return $\Sigma_{-\frac{1}{2}}$
			\EndFunction
		\end{algorithmic}
	\end{samepage}
\end{onehalfspace}
\vspace{2 mm}

\section{Simulations}

We shall now show some simulation results of adaptive detectors using various estimators introduced in this article, using the detection tests introduced in \cite{CVC_BT}\cite{scharf1994matched}.\\

The simulation results are shown in a single channel scenario of dimension $d=8$, as a function of the SiNR of the target signal.\\
the background noise is generated as a white gaussian noise of unit variance.\\
The target signal is generated as a complex centered circular 1-dimensional gaussian signal aligned with the test signal $s$, whose variance $\sigma$ is such that:

\begin{displaymath}
10 \log_{10}(\sigma) = \text{SiNR}
\end{displaymath}

The detection thresholds are defined to have a false alarm rate of $10^{-4}$; they are learned on clean training sets, that is that they contain no outliers.\\

In this scenario, we have an infinite stream of learning samples 
$(x_{n,m})_{1\leq n \leq N, m \in \mathbb{Z}}$. These learning samples are used to learn the covariance/correlation matrix of the signal sequentially, using the following pattern:

\begin{displaymath}
\Sigma_m = \text{reg}(\Sigma_{m-1},(x_{n,m})_{1\leq n \leq N})
\end{displaymath}

This covariance/correlation matrix is then used in an adaptive detection scheme.\\

The performances are considered in the asymptotic regime (meaning that $\Sigma_m$ follows the same distribution as $\Sigma_{m-1}$ for any $m \in \mathbb{Z}$). The performances are shown for $N = 11$, $p = \frac{3}{4}$ for partial estimators, with various values of $\alpha$; the results are also compared to the non regularized versions of the estimators for $N = 22$ \cite{CVC_BT}\cite{CVC_partial}.

\begin{figure}[H]
	\includegraphics[width=\linewidth]{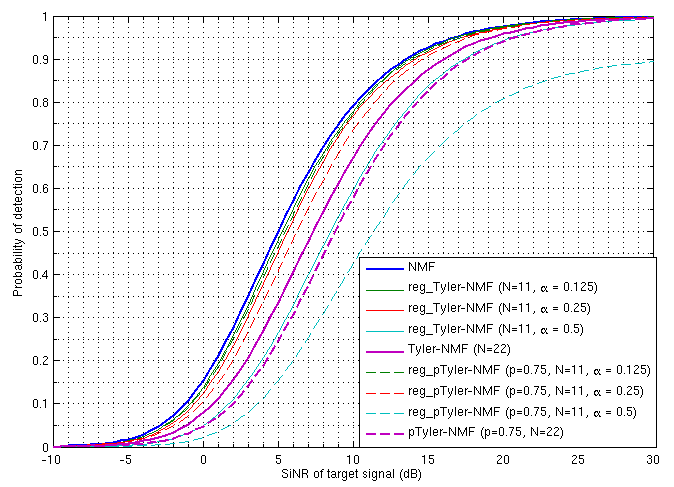}
	\caption{Detection capabilities of Tyler-NMF and pTyler-NMF}
	\label{fig_reg_Tyler_pTyler}
\end{figure}

\begin{figure}[H]
	\includegraphics[width=\linewidth]{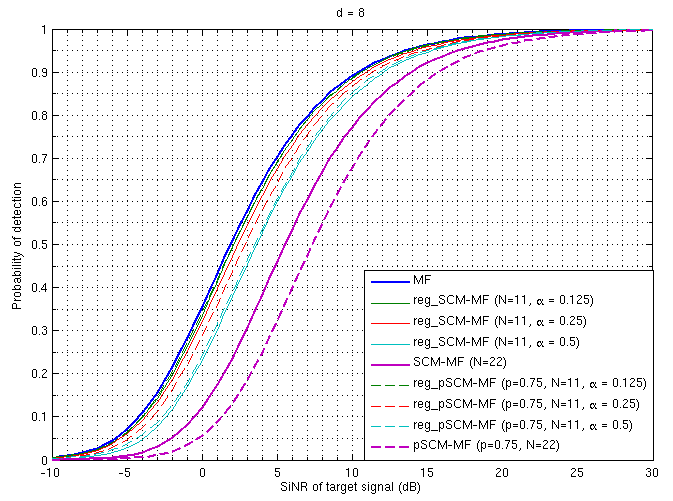}
	\caption{Detection capabilities of reg\_SCM-MF and reg\_pSCM-MF}
	\label{fig_reg_scm}
\end{figure}

\begin{figure}[H]
	\includegraphics[width=\linewidth]{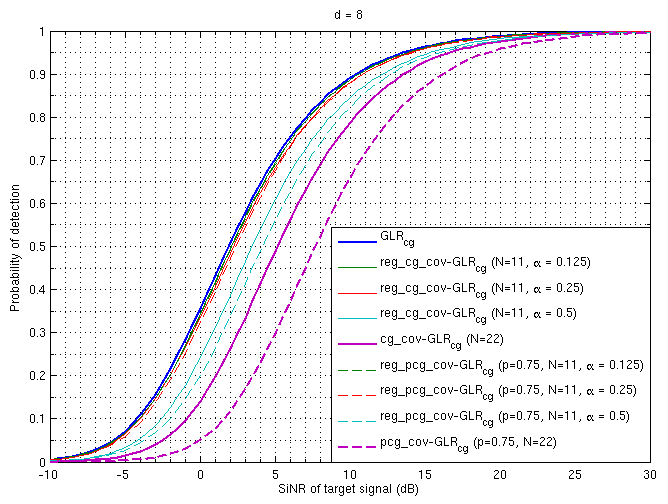}
	\label{fig_reg_cg_pcg}
\end{figure}

As can be seen on figures \ref{fig_reg_Tyler_pTyler}, \ref{fig_reg_scm} and \ref{fig_reg_cg_pcg}, the use of a regularized estimate in this case allows to greatly improve the detection capabilities of adaptive tests. Indeed the performances quickly approach the asymptotic limit in which the covariance matrix is perfectly estimated when a slower integration factor is chosen.\\
Moreover the degradation of performances of partial estimators can be almost completely negated for a slow enough integration ($\alpha \leq \frac{1}{4}$), thus resulting in estimators which offer performances close to non-adaptive detectors, while retaining the robustness to outliers offered by partial estimators.

\section{Conclusion}

We have presented a theoretical background for a likelihood based approach to regularized estimations. This background allows to define penalties used in penalized maximum likelihood schemes in a way which is consistent with the maximum likelihood principle.\\
This maximum likelihood principle for regularized estimation is then further extended to allow for robust estimators in the presence of outliers, by using the partial procedure defined in \cite{CVC_partial}.\\
Finally the potential gain in terms of detection performances in adaptive detection schemes is shown on some simulations.

\appendices
\section{On the relative entropy of angular models}
\label{sec:KL_angular}
The relative entropy of two angular models of $\mathbb{R}^d$, characterized by their correlation matrix $R_0$ and $R$ (with the convention that $|R_0|=|R|=1$), is given by:
\begin{flalign*}
H(R_0|R) = &&
\end{flalign*}
\begin{equation}
\frac{d-1}{2}\int_{x \in \mathcal{S}_{d-1}} \log\left(\frac{x^\dagger R^{-1} x}{x^\dagger {R_0}^{-1} x}\right)ds_{d-1}\left(\frac{L(R_0)^{-1}x}{\sqrt{x^\dagger {R_0}^{-1} x}}\right)
\end{equation}

Although the exact computation of this integral represents a challenge, what is of interest to us is its differential with respect to the matrix $R$:
\begin{flalign*}
dH(R_0|R) = &&
\end{flalign*}
\begin{equation}
\frac{d-1}{2}\text{tr}\left(\int\limits_{x\in S_{d-1}}  \frac{xx^\dagger}{x^\dagger R^{-1}x}ds_{d-1}\left(\frac{L(R_0)^{-1}x}{\sqrt{x^\dagger {R_0}^{-1} x}}\right)d(R^{-1})\right)
\end{equation}

Thus we have to compute:

\begin{equation}
\Sigma = \int_{x\in S_{d-1}} \frac{xx^\dagger}{x^\dagger R^{-1}x}ds_{d-1}\left(\frac{L(R_0)^{-1}x}{\sqrt{x^\dagger {R_0}^{-1} x}}\right)
\end{equation}

By making a change of variable $y = R^{-\frac{1}{2}}x$, this integral can be recast as:

\begin{equation}
\Sigma = R^\frac{1}{2} \left(\int_{y \in S_{d-1}} \frac{yy^\dagger}{y^\dagger y} ds_{d-1}\left(\frac{L(R^{-\frac{1}{2}}R_0 R^{-\frac{1}{2}})^{-1}y}{\sqrt{y^\dagger R^{-\frac{1}{2}}R_0 R^{-\frac{1}{2}}y}}\right)\right)R^\frac{1}{2}
\end{equation}

Let us note for any symmetric positive matrix $\Sigma$:

\begin{equation}
\text{anscm}(\Sigma) = \int_{y \in S_{d-1}} \frac{yy^\dagger}{y^\dagger y}ds_{d-1}\left(\frac{L(\Sigma)^{-1}y}{\sqrt{y^\dagger \Sigma^{-1}y}}\right)
\end{equation}

This is the integral of $\frac{yy^\dagger}{y^\dagger y}$ on the unit ellipse defined by the correlation matrix $\Sigma$, and therefore corresponds to the  expected value of the normalized sample covariance estimator (NSCM) for an elliptical distribution of correlation matrix proportional to $\Sigma$. The structure of this matrix has been studied thoroughly in \cite{bausson2007first} for complex circular distributions, in the case in which all eigenvalues are distinct: 

\begin{equation}
\text{anscm}(\Sigma) =  U\text{diag}(\mu)U^\dagger
\end{equation}

with: 

\begin{itemize}
	\item $\Sigma = U \text{diag}(\lambda) U^\dagger$ with $U U^\dagger = I$ being the eigen decomposition of $\Sigma$
	\item $\mu = Ec$ with:
	\begin{displaymath}
	\left\{\begin{array}{l}
	E = \left(\left\{\begin{array}{lc}
	\frac{\log(\frac{\lambda_j}{\lambda_i})}{\frac{\lambda_j}{\lambda_i}-1}-\frac{\lambda_i}{\lambda_j} & \mbox{if } i \neq j\\
	0 & \mbox{if } i = j
	\end{array}\right.\right)_{1 \leq i,j \leq d}\\
	c = \left(\prod\limits_{k \neq j} \frac{1}{1-\frac{\lambda_k}{\lambda_j}}\right)_{1 \leq j \leq d}
	\end{array}\right.
	\end{displaymath}
\end{itemize}

Note that in a practical implementation, the limitation of having distinct eigenvalues is not restrictive, as it is almost surely the case. However one should take steps to insure a decent numerical stability. We have on our end opted for enforcing a minimum spacing between eigenvalues determined by a numerical tolerance, as given by the following algorithm:\\

\begin{onehalfspace}
	\begin{samepage}
		\begin{algorithmic}[1]
			\Function{Respace}{$\left(\lambda_k\right)_{1 \leq k \leq d},\epsilon$}
			\State $o \gets \text{argsort}_\uparrow\left(\left(\lambda_k\right)_{1 \leq k \leq d}\right)$\vspace{1.5 mm}
			\State $\left(\mu_k\right)_{1\leq k\leq d} \gets \left(\lambda_k\right)_{1 \leq k \leq d}$\vspace{1.5 mm}
			\For{$k$ \textbf{from} $1$ \textbf{to} $d-1$}
			\If{$\frac{\mu_{o(k+1)}}{\mu_{o(k)}}-1<\epsilon$}\vspace{1.5 mm}
			\For{$m$ \textbf{from} $k+1$ \textbf{to} $d$}
			\State $\mu_{o(m)} \gets (1+\epsilon)\mu_{o(m-1)}-\mu_{o(m)}$\vspace{1.5 mm}
			\EndFor
			\EndIf
			\EndFor
			\State \Return $\left(\mu_k\right)_{1\leq k\leq d}$
			\EndFunction
		\end{algorithmic}
	\end{samepage}
\end{onehalfspace}
\vspace{2 mm}

Unfortunately we do not have a similar result in the general case of distributions in $\mathbb{R}^d$, in which case one has to resort to numerical integration methods. Indeed the matrix can be expressed as $\text{anscm}\left(U\text{diag}(\lambda)U^\dagger\right) = U\text{diag}(\mu)U^\dagger$, with:

\begin{flalign*}
\mu_k = &&
\end{flalign*}
\begin{equation}
\frac{2}{\sqrt{\pi}}\frac{\Gamma\left(\frac{d}{2}\right)}{\Gamma\left(\frac{d-1}{2}\right)}\iint\limits_{\substack{
		t\in[0;1]\\
		x\in \mathcal{S}_{d-2}
	}
}
\frac{{t^2(1-t^2)}^{\frac{d}{2}-1}}{t^2+\left(1-t^2\right)\sum\limits_{m \neq k} \frac{\lambda_m}{\lambda_k}{x_k}^2}ds_{d-2}(x)dt
\end{equation}

There is no known analytical expression of these eigenvalues to the authors' knowledge.



\ifCLASSOPTIONcaptionsoff
  \newpage
\fi



%

\bibliographystyle{plain}
\bibliography{CES_biblio}



%

\begin{IEEEbiographynophoto}{Christophe Culan}
	was born in Villeneuve-St-Georges, France, on November $\text{1}^\text{st}$, 1988. He received jointly the engineering degree from Ecole Centrale de Paris (ECP), France and a master's degree in Physico-informatics from Keio University, Japan, in 2013.\\
	He was a researcher in applied physics in Itoh laboratory from 2013 to 2014, specialized in quantum information and quantum computing, and has contributed to several publications related to these subjects.\\
	He currently holds a position as a research engineer in Thales Air Systems, Limours, France, in the Advanced Radar Concepts division. His current research interests include statistical signal and data processing, robust statistics, machine learning and information geometry.
\end{IEEEbiographynophoto}

\begin{IEEEbiographynophoto}{Claude Adnet}
	was born in Aÿ, France, in 1961. He received the DEA degree in signal processing and Phd degree in 1988 and 1991 respectively, from the Institut National Polytechnique de Grenoble (INPG), Grenoble France.
	Since then, he has been working for THALES Group, where he is now Senior Scientist. His research interests include radar signal  processing and radar data processing.
\end{IEEEbiographynophoto}






\end{document}